\documentclass[12pt]{amsart} % use larger type; default would be 10pt

\usepackage{amssymb, latexsym}

\newtheorem{Theorem}{Theorem}

\newtheorem{Lemma}[Theorem]{Lemma}
\newtheorem{Proposition}[Theorem]{Proposition}

\usepackage%[backref]
{hyperref}
\usepackage{graphicx}

\begin{document}

\title[Lerch zeta-function for equal parameters] {Zeros of the Lerch zeta-function and of its derivative for equal parameters}

\author{ Ram\={u}nas Garunk\v{s}tis}
\address{Ram\={u}nas Garunk\v{s}tis \\
Department of Mathematics and Informatics, Vilnius University \\
Naugarduko 24, 03225 Vilnius, Lithuania}
\thanks{The first author is  supported
by grant No. MIP-049/2014 from the Research 
Council of Lithuania.}
\email{ramunas.garunkstis@mif.vu.lt}
\urladdr{www.mif.vu.lt/~garunkstis}

\author{ Rokas Tamo\v si\=unas}
\address{ Rokas Tamo\v si\=unas \\
Department of Mathematics and Informatics, Vilnius University \\
Naugarduko 24, 03225 Vilnius, Lithuania}
\email{trokas@gmail.com}

\subjclass[2010]{Primary: 11M35; Secondary: 11M26}

\keywords{Lerch zeta-function; nontrivial zeros; Speiser's equivalent for the Riemann hypothesis}

\begin{abstract} 
A. Speiser proved that the Riemann hypothesis is equivalent to the absence of non-real zeros of the derivative of the Riemann zeta-function left of the critical line. His result has been extended by N. Levinson and H.L. Montgomery to the statement that the
Riemann zeta-function and its derivative have approximately the same number of non-real zeros  left of the critical line. 
We obtain the  Levinson-Montgomery type result for the Lerch zeta-function with equal parameters. For the Lerch zeta-function, the analogue of the Riemann hypothesis is usually not true and its zeros usually are distributed asymmetrically with respect to the critical line. However, for equal parameters, the symmetry of the zeros is almost restored.
\end{abstract}

\maketitle

\section{Introduction}

Let $s=\sigma+it$. Denote by $\{\lambda\}$ the fractional part of a real number $\lambda$. In this paper $T$ always tends to plus infinity. In all theorems and lemmas, the numbers $\lambda$ and $\alpha$ are fixed constants.

For $0<\lambda, \alpha\leq 1$, the Lerch zeta-function is given by
\begin{align*}
L(\lambda,\alpha,s)=\sum_{m=0}^\infty\frac{e^{2\pi i\lambda m}}{(m+\alpha)^s} \qquad(\sigma>1).
\end{align*}
This function has analytic continuation to the whole complex plane except for a possible simple pole at $s=1$ (Lerch \cite{Lerch1887}, Laurin\v cikas and Garunk\v stis  \cite{gl}).

Let $\zeta(s)$ and $L(s,\chi)$ denote the Riemann zeta-function and the Dirichlet $L$-function accordingly. We have that $L(1,1,s)=\zeta(s)$ and $L(1/2,1/2,s)=2^s L(s,\chi)$, where $\chi$ is a Dirichlet character$\mod 4$
with $\chi(3)=-1$.  For these two cases, certain versions of the Riemann hypothesis (RH) can be formulated. Similar cases are $ L(1,1/2,s)=(2^s-1)\zeta(s)$ and  $L(1/2,1,s)=(1-2^{1-s})\zeta(s)$. 
%\begin{align*}
%&L(1,1,s)=\zeta(s),\quad L(1,1/2,s)=(2^s-1)\zeta(s),\\
%&L(1/2,1,s)=(1-2^{1-s})\zeta(s),\quad\mathrm{and}\quad
%L(1/2,1/2,s)=2^s L(s,\chi), \nonumber
%\end{align*}
% For these four cases, certain versions of the Riemann hypothesis can be formulated.  
For all the other cases,  it is expected that the real parts of zeros of the Lerch zeta-function form a dense subset of the interval $(1/2,1)$. This is proved for any $\lambda$ and transcendental $\alpha$ (\cite[Theorem 4.7 in Chapter 8]{gl}). 

Speiser \cite{Speiser1934} showed that the Riemann hypothesis (RH) is
equivalent to the absence of non-real zeros of the derivative of the
Riemann zeta-function left of the critical line. Later on, Levinson and Montgomery \cite{Levinson1974} proved the quantitative version of the
Speiser's result, namely, that the Riemann zeta-function and its
derivative have approximately the same number of zeros left of the
critical line. This result was extended to Dirichlet $L$-functions with
primitive Dirichlet characters (Y{\i}ld{\i}r{\i}m \cite{yildirim96}), to the
Selberg class (\v Sle\v zevi\v cien\. e \cite{rasa}), to the Selberg
zeta-function on a compact Riemann surface (Luo \cite{Luo2005},  Garunk\v stis \cite{Garunkstis2008}). See also Minamide \cite{Minamide2009}, \cite{Minamide2010}, \cite{Minamide2013},  Jorgenson and Smailovi\' c \cite{js}. In  
these cases an
analog of the RH is expected or, as in the case of the Selberg
zeta-function  on a compact Riemann surface, it is known to be true. In Garunk\v stis and \v Sim\.enas \cite{gs}, the Speiser equivalent was investigated for the extended Selberg class. Zeta-functions of this class satisfy  functional equations of classical type, thus the nontrivial zeros are  distributed symmetrically with respect of the critical line. Moreover, this class contains zeta-functions for which the analog of RH is not true. 

Here we consider the relation between zeros of the Lerch zeta-function and its derivative  when parameters are equal. 
In our paper \cite{gt} we    showed  that the nontrivial zeros of $L(\lambda, \lambda, s)$  either lie extremely close to the critical line $\sigma = 1/2$
or are distributed almost symmetrically with respect to the critical line. Detailed calculations however suggest that this symmetry is not strict if $0<\lambda<1$ and $\lambda\ne1/2$.

 For the Lerch zeta-function the following relation, usually called the functional equation, is true.
\begin{align}\label{Lerchfunc}
L(\lambda,\alpha,1-s)=&(2\pi)^{-s}\Gamma(s)\biggr(
e^{\pi i\frac{s}{2}-2\pi i\alpha\lambda}L(1-\alpha,\lambda,s)
\\&
+e^{-2\pi i\frac{s}{4}+2\pi i \alpha(1-\{\lambda\})}L(\alpha,1-\{\lambda\},s)\biggr).\nonumber
\end{align}
Various proofs of this functional equation can be found in Lerch \cite{Lerch1887}, Apostol \cite{Apostol1951}, Oberhettinger \cite{Oberhettinger1956}, Mikol\'as \cite{Mikolas1971}, Berndt \cite{Berndt1972}, see also Lagarias and Li \cite{ll1}, \cite{ll2}.  The almost symmetrical distribution of zeros in the case of equal parameters  is related to the functional equation 
 \eqref{Lerchfunc}, which for $\lambda=\alpha$ can be rewritten as 
 \begin{align*}%\label{almostsymetry}
\overline{L(\lambda,\lambda,1-\overline{s})}=&(2\pi)^{-s}\Gamma(s)e^{-\pi i\frac{s}{2}+2\pi i \lambda^2}L(\lambda,\lambda,s)\nonumber
\\&
+(2\pi)^{-s}\Gamma(s)
e^{\pi i\frac{s}{2}-2\pi i(1-\lambda)\lambda}L(1-\lambda,1-\{\lambda\},s)
\\
=&G(s)L(\lambda,\lambda,s)+P(s),\nonumber
\end{align*}
where, for any vertical strip, $|P(s)|< t^Be^{-\pi t}$ and $|G(s)|\ge t^C$, $B, C>0$ (see \cite{gt}).

Next, we recall several facts about the Lerch zeta-function. Later, we  formulate the obtained result.

 Let $l$ be a straight line in the complex plane ${\Bbb C}$,
and denote by $\varrho(s,l)$ the distance of $s$ from $l$.
Define, for $\delta>0$,
\begin{eqnarray*}
L_\delta(l)=\big\{s\in{\Bbb C}:\;\varrho(s,l)<\delta\big\}.
\end{eqnarray*}
In Garunk\v stis and Laurin\v cikas \cite{gar}, Garunk\v stis and Steuding \cite{Garunkstis2002}, for $0<\lambda<1$ and $\lambda\ne1/2$,   it is proved that 
$L(\lambda,\alpha,s)\ne0$ if $\sigma<-1$ and
\begin{eqnarray*}
s\not\in L_{\log4\over\pi}\bigg(\sigma =\frac{\pi
t}{\log{\frac{1-\lambda}{\lambda}}}+1\bigg)\,.
\end{eqnarray*}
  For $\lambda=1/2,1$, from Spira \cite{Spira} and \cite{gar} we see that  $L(\lambda,\alpha,s)\ne0$ if $\sigma<-1$ and $|t|\ge1$.  Moreover, in  \cite{gar} it is showed that $L(\lambda,\alpha, s)\ne0$ if $\sigma\ge1+\alpha$.  We say that a zero of $L(\lambda,\alpha, s)$ is {\it nontrivial} if it lies in the strip $-1\le\sigma<1+\alpha$ and we denote a nontrivial zero by $\rho=\beta+i\gamma$.

Denote by $N(\lambda,\alpha,T)$  the number of 
nontrivial zeros of the function $L(\lambda,\alpha,s)$ in the region 
$0<t<T$. 
For $0<\lambda,\alpha\leq 1$,  we have (\cite{gar}, Garunk\v stis and Steuding \cite{gs2002})
\begin{align}\label{zeronumber}
N(\lambda,\alpha,T)=\frac{T}{2\pi}\log 
\frac{T}{2\pi e\alpha\lambda}+O(\log T).
\end{align}
For recent results on the value-distribution of the Lerch zeta-function see  Mishou \cite{Mishou14},  Lee, Nakamura, and Pa\'nkowski \cite{Lee17}.

Define
\begin{align*}
L'(\lambda,\lambda,s)={\partial\over\partial s}L(\lambda,\lambda,s). 
\end{align*}
We collect several facts about the zero distribution of $L'(\lambda,\lambda,s)$.  From the expression of $L'(\lambda,\lambda,s)$ by the Dirichlet series, we have that there is $\sigma_1\ge1$ such that
 $L'(\lambda,\lambda,s)\ne0$ 
if $\sigma>\sigma_1$. By Lemma \ref{lefthandside} below and the zero free region of $L(\lambda,\lambda,s)$ we see that, for any $\sigma<-1$, there is a constant $t_0=t_0(\sigma)$ such that $$L'(\lambda,\lambda,s)\ne0$$ if $t\ge t_0$.
 We say that a zero of $L'(\lambda,\lambda, s)$ is {\it nontrivial} if it lies in the strip $-1\le\sigma\le\sigma_1$.  Let $N'(\lambda,\lambda,T)$ denote the number of 
nontrivial zeros of  $L'(\lambda,\lambda,s)$ in the region
$0<t<T$. We have (\cite[Notes to Chapter 8]{gl})
\begin{align*}
N'(\lambda,\lambda, T)={T\over {2\pi}}\log{T\over 2\pi e
([\lambda]+\lambda)\lambda}+o(T).
\end{align*}
Therefore $N(\lambda,\lambda, T)-N'(\lambda,\lambda, T)=o(T)$ whenever $\lambda\neq 1$.

To state our main result, we need several notations. By the formula (\ref{zeronumber}),  there is a constant $D=D(T_0)$ such that, for $T>T_0$,
\begin{align}\label{D}
\left|N(\lambda,\lambda,T)-\frac{T}{2\pi}\log 
\frac{T}{2\pi e\lambda^2}\right| \le D\log T.
\end{align}
Let $\psi :[0,1]\to\mathbb C$ always denote a simple piecewise smooth curve with the initial point on the line $t=T$ and the terminal point on the line $t=T+U$. Moreover, let $-2<\Re \psi(\tau)\le1/2$ and $T<\Im \psi(\tau)<T+U$, where $\tau\in(0,1)$. 
Let $M(T, U, \psi)$ (resp. $M'(T, U, \psi)$) be the number of nontrivial zeros of $L(\lambda,\lambda, s)$ (resp. $L'(\lambda,\lambda, s)$) inside (but not on the border of) the area restricted by segments $[1/2+i(T+U), -2+i(T+U)]$, $[-2+i(T+U), -2+iT]$, $[-2+iT, 1/2+iT]$, and the curve $\psi$.
\begin{Theorem}\label{dertest}
Let $0<\lambda\le1$ and  $T>0$.    
Assume that, for some $T_0$ and $0<\varepsilon<1$,
\begin{align}\label{Dcondition}
D(T_0)<\frac{\varepsilon}{\log 2}. 
\end{align} 
Then, for sufficiently large $T$ and $0<U\le T$, there is a positive constant $A$ and a curve $\psi :[0,1]\to\mathbb C$ such that, 
\begin{align*}
1/2-\exp(-AT^{1-\varepsilon}/\log T)\le \Re\psi(\tau)\le1/2\quad(\tau\in[0,1]),
\end{align*}
and
$$M(T, U, \psi)=M'(T, U,\psi)+O(\log T).$$
\end{Theorem}

We discuss the condition (\ref{Dcondition}). For the Riemann zeta-function ($=L(1,1,s)$) it is known that $D<0.12<1/\log 2=1.44\dots$ (Trudgian \cite{Trudgian2014}). If $\lambda=1/2$, then $D<0.16$ (Trudgian \cite{Trudgian2015}).  Moreover, for the Riemann zeta-function the Lindel\"of hypothesis implies that the constant $D$ can be chosen as small as we please (Titchmarsh \cite[Theorem 13.6(A)]{Titchmarsh1986}).  We expect the Lindel\"of type hypothesis also for the Lerch zeta-function (\cite{Garunkstis2002},  \cite{Garunkstis2005}). Similarly as in the case of the Riemann zeta-function (\cite[Sections 13.6 and 13.7]{Titchmarsh1986}), it is possible to modify the proof of Theorem 3.2 in \cite[Chapter 8]{gl} and to show that the  Lindel\"of type hypothesis for $L(\lambda,\alpha,s)$ implies that, for any $0<\lambda<1$, the constant $D$  can be chosen as small as we please. 

In the next section we present the computer computations related to  Theorem \ref{dertest}. Section  \ref{proofTh4} contains the proof  of Theorem   \ref{dertest}. In the last section we discuss the curve $\psi$ from Theorem \ref{dertest}.

\section{Computations}

Here we draw several graphs illustrating Theorem \ref{dertest}.
We consider the  trajectories of the zeros of $L(\lambda, \lambda, s)$ and of its derivative. 

Suppose that $\rho=\rho(\lambda_0)$ is a zero of multiplicity $m>0$ of
$L(\lambda_0, \lambda_0, s)$. From 
the expression of the Lerch zeta-function by the Dirichlet series and from the functional equation (\ref{Lerchfunc}), it follows that, for any $s$, the function $f(\lambda)=L(\lambda, \lambda, s)$ is continuous in $\lambda\in(0,1)$. By Rouch\'
e's theorem, we have that for every sufficiently small open disc $D$
with center at $\rho$ in which the function $L(\lambda_0, \lambda_0,s)$ has no other
zeros except for $\rho$, there exists $\delta=\delta(D)>0$ such that
each function $L (\lambda, \lambda, s)$, where $\lambda\in (\lambda_0-\delta,
\lambda_0+\delta)$, has exactly $m$ zeros (counted with multiplicities)
in the disc $D$ (c.f. Theorem 1 in Balanzario and S\'{a}nchez-Ortiz
\cite{Balanzario2007} and Lemma 4.1 in Dubickas, Garunk\v stis,
J. Steuding and R. Steuding \cite{DGSS}). If zero $\rho$ is of
multiplicity $m=1$, then there exists a neighborhood of $\lambda_0$ and
some function $\rho=\rho(\lambda)$, which is continuous at $\lambda_0$ and,
in addition, satisfies the relation $L(\lambda, \lambda, \rho(\lambda))=0$. This way,
we can speak about the continuous trajectory $\rho(\lambda)$ of a zero. The trajectories of the zeros of the
derivative
$L'(\lambda,\lambda,s)$ are understood in a similar way. 

\begin{figure}[h]
\includegraphics[width=0.32\textwidth]{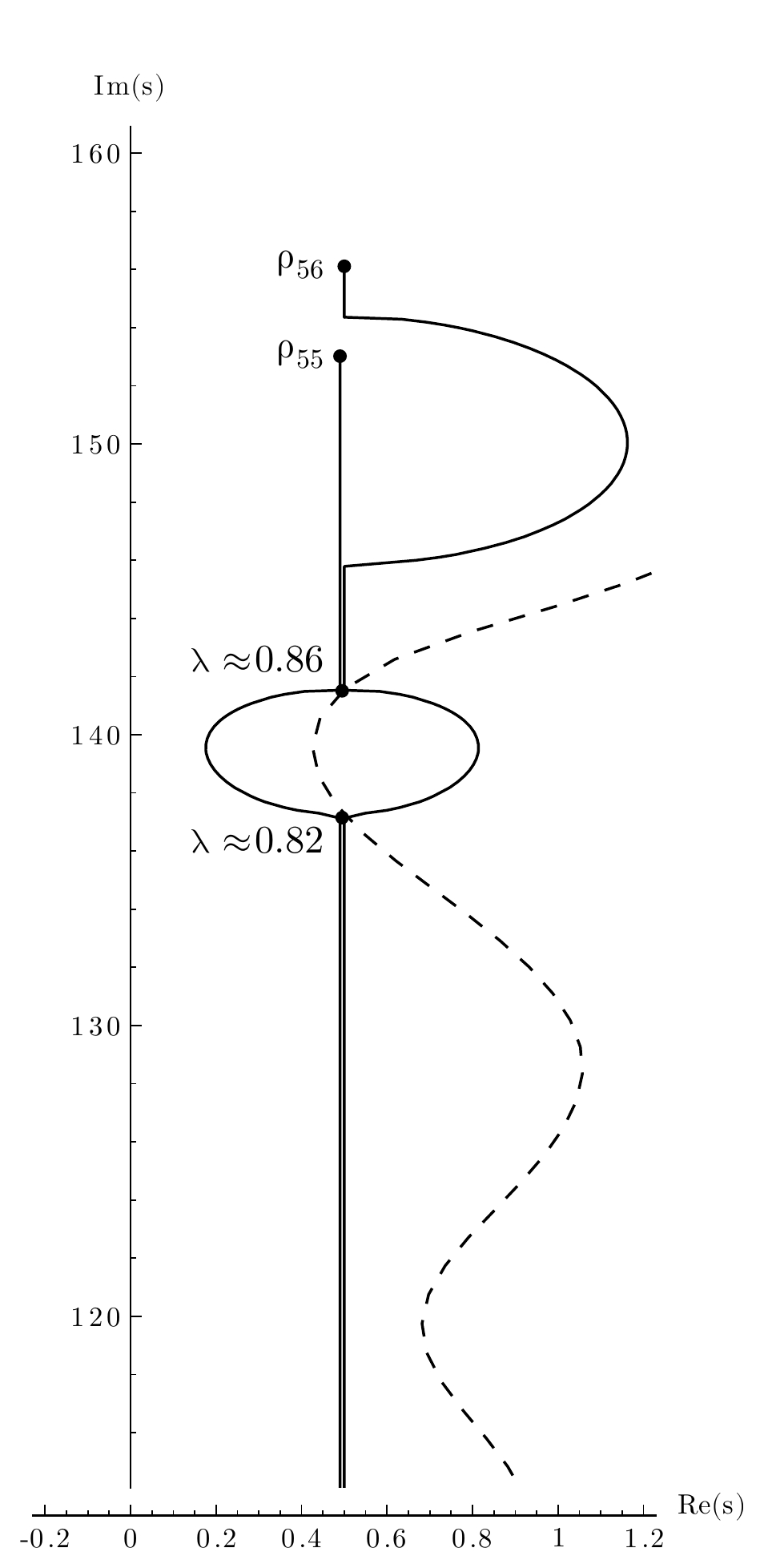}
\hfill
\includegraphics[width=0.32\textwidth]{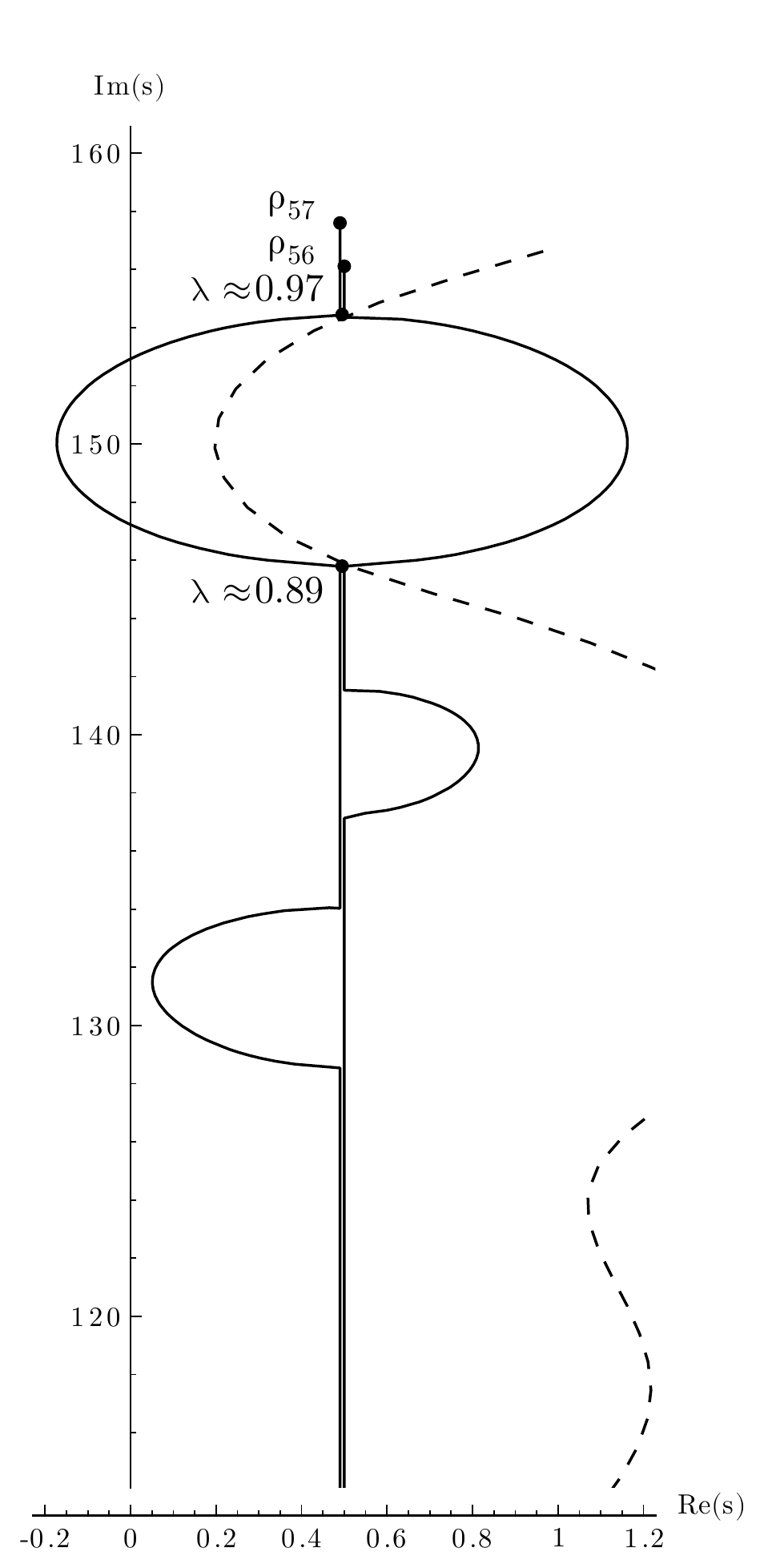}
\hfill
\includegraphics[width=0.32\textwidth]{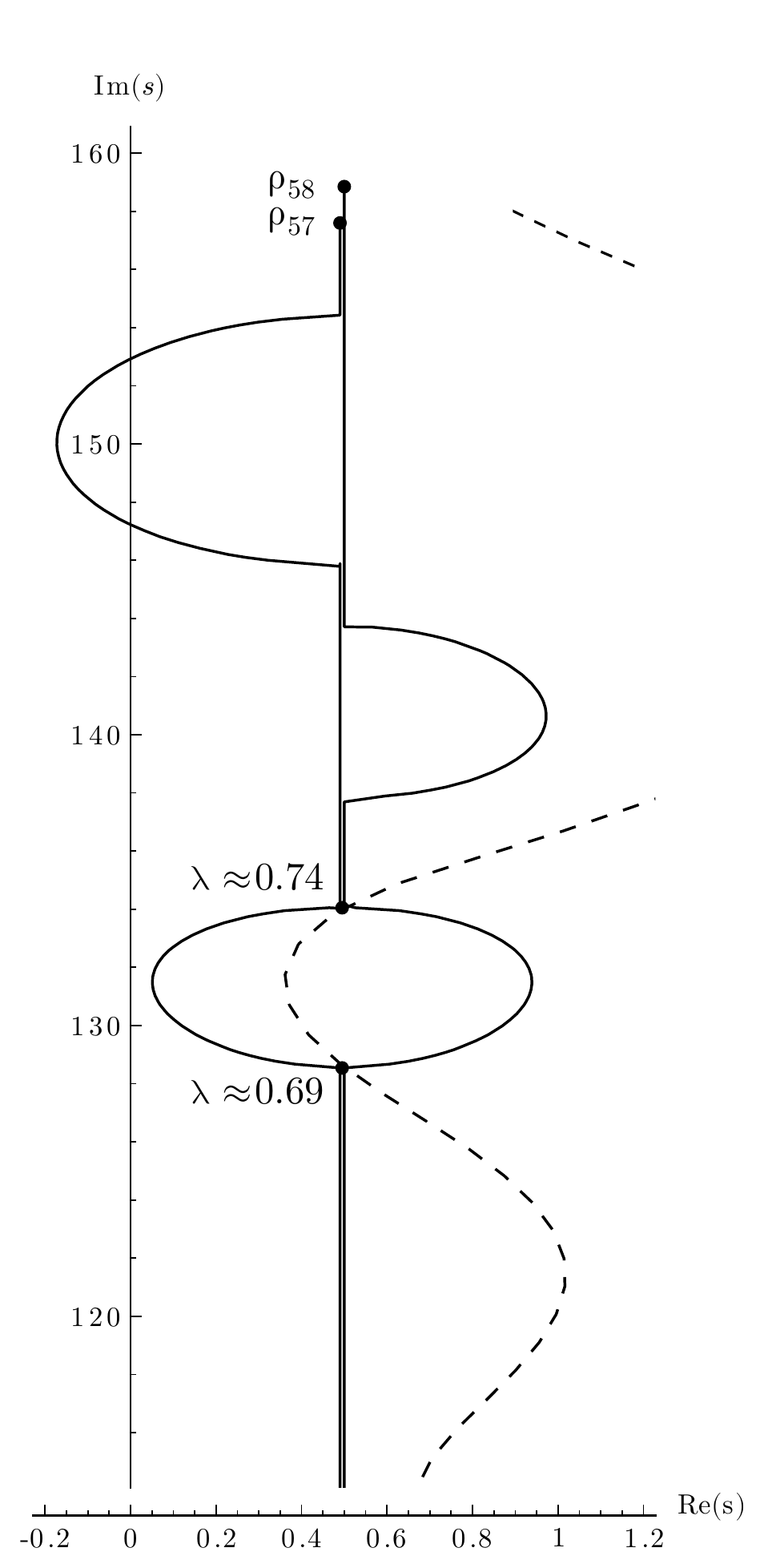}
\\
\caption{
Parametric graphics of the several  trajectories of the zeros with parameter $1/2\le\lambda\le1$. Solid  and dashed trajectories are the trajectories of the zeros of   $L\left(\lambda,\lambda, s\right)$ and $L'\left(\lambda,\lambda, s\right)$  respectively. By $\rho_{55},\dots,\rho_{58}$ we denote  55th,\dots, 58th zeros of $\zeta(s)=L(1,1,s)$. The trajectories of the zeros of the derivative correspond to the zeros   $1.27.. + 152.61..i$ (left), $0.97.. + 156.63..i$ (middle), and $0.86.. + 158.28..i$ (right)  of  $\zeta'(s)=L'(1,1,s)$.
}\label{fig2}
\end{figure}

In Figure \ref{fig2},  we see parametric
plots of the trajectories of the zeros of $L(\lambda,\lambda,s)$ and its
derivative, solid and dotted lines respectively. We see that the
trajectory of the derivative crosses the critical line in accordance
with Theorem~\ref{dertest}. 

Figure \ref{fig2} can be compared to Figures 1 and 2 in Garunk\v stis and \v Sim\. enas \cite{gs}, where the trajectories of the zeros of the linear combination $f(s,\tau)$
of Dirichlet $L$-functions and of its derivative $f'_s(s,\tau)$ is calculated. Here, in Figure \ref{fig2}, we see that  the trajectories of the  zeros of $L(\lambda,\lambda,s)$ approach the (almost) meeting point (after which the trajectories leave a neighborhood  of the critical line) from the same direction, while   the trajectories of the zeros of $f(s,\tau)$ approach the  meeting point from the opposite directions. Note that the nontrivial zeros of this  linear combination of Dirichlet $L$-functions  are distributed strictly symmetrically with respect of the critical line. Because of this fact, the meeting point in Figures 1 and 2 in  \cite{gs} is always a double zero of $f(s,\tau)$.

To find  the trajectories of the zeros $\rho(\lambda)$ and $q(\lambda)$, 
$0<\lambda\le1$, such that
$$L(\lambda, \lambda, \rho(\lambda))=0\quad\text{and}\quad L'(\lambda, \lambda, q(\lambda))=0,$$
we solve the differential equations numerically
\begin{align*}
  \frac{\partial {\rho}(\lambda )}{\partial \lambda }=
  {{-\frac{\frac{\partial \ell (\lambda, \rho)}{\partial \lambda
        }}{\frac{\partial \ell  (\lambda, \rho)}{\partial
          \rho}}}}\quad\text{and}\quad\frac{\partial {q} (\lambda)}{\partial \lambda }= {{-\frac{\frac{\partial^2 \ell(\lambda, q)}{\partial q\partial \lambda }}{\frac{\partial^2 \ell
          (\lambda, q)}{\partial q^2}}}},
\end{align*}
where $\ell(\lambda, s)=L(\lambda,\lambda,s)$.
As the initial conditions, some zeros of $L(1,1,s)= \zeta (s)$ and $L'(\lambda,\lambda,s)$, $\lambda=0.86.., 0.97.., 0.74..$ are used.

Computations were performed by {\it Wolfram Mathematica} using commands {\it LerchPhi}, {\it NDSolve},  also numerically calculating contour integrals 
$$\oint\frac{L'(\lambda,\lambda,s)}{L(\lambda,\lambda,s)}ds\quad \text{and}\quad \oint s\frac{L'(\lambda,\lambda,s)}{L(\lambda,\lambda,s)}ds.$$ 

Computations were validated with the help of Python with {\it mpmath}\footnote{Fredrik Johansson and others. mpmath: a Python library for arbitrary-precision floating-point arithmetic (version 0.18), December 2013. http://mpmath.org/.} package.  We used the following expression of the Lerch zeta-function for rational parameters  %(see \eqref{eq:expr}) 
\begin{align*}
L\left(\frac{b}{d},\frac{b}{d}, s\right) & =\sum_{k=0}^{d-1}\sum_{m=0}^{\infty}\frac{\exp\left(2\pi i\frac{b}{d}\left(dm+k\right)\right)}{\left(dm+k+\frac{b}{d}\right)^{s}}\\
& =d^{-s}\sum_{k=0}^{d-1}\exp\left(2\pi i\frac{b}{d}k\right)\zeta\left(s,\frac{kd+b}{d^{2}}\right),
\end{align*}
where $\zeta(s,\alpha)$, $0<\alpha\le1$, is the Hurwitz zeta-function.  The function $\zeta(s,\alpha)$ is
implemented by  the command {\it zeta}.  Zero locations were calculated using {\it findroot} with Muller\textquoteright s method. 

In this paper, all computer computations  should be regarded as heuristic because their accuracy was not controlled explicitly.

\section{Proof of Theorem \ref{dertest}}\label{proofTh4}

The structure of the proof is similar to the proof of the formula (10.28.2) in  Section 10.28 of Titchmarsh \cite{Titchmarsh1986}, see also the original proof in Levinson and Montgomery \cite{Levinson1974}. The main difference is Proposition \ref{crho} below.

The following bound from below for $L(\lambda, \lambda, s)$ when $s$ is close to a zero will be useful. 
\begin{Lemma}\label{1/L}
Let $0<\lambda, \alpha\le1$. Let $\sigma_2\in\mathbb R$ and $\Re s\ge\sigma_2$. Let $L(\lambda,\alpha,s)\ne0$ and  $d$ be the distance from $s$ to the nearest zero of $L(\lambda,\alpha,s)$. Then, for $t\ge2$,
\begin{align*}
\frac1{|L(\lambda,\alpha, s)|}<\exp(C(|\log d|+1)\log t),
\end{align*}
where $C=C(\lambda, \alpha, \sigma_2)$ is a positive constant.
\end{Lemma}
\proof This is  Proposition 1 in our paper \cite{gt}.

 \endproof

\begin{Lemma}\label{lefthandside}
Let $0<\lambda\le1$, $T<t<T+U$, and $0<U\le T$. If  $\sigma_2<-1$ then
\begin{align}\label{left}
\Re\frac{L'}{L}(\lambda,\lambda,\sigma_2+it)=-\log t+O_{\sigma_2}(1)\qquad(T\to\infty).
\end{align}

Moreover, assume that $A>0$ is such that $4AC<\pi$, where $C=C(\lambda,\lambda, 1/2)$ is a constant from Lemma \ref{1/L}. If the distance from $1/2+it$ to the nearest zero of  $L(\lambda,\lambda, s)$ is greater than  $\exp\left(-A T/\log T\right)$, then
\begin{align}\label{logr11/2}
\Re \frac{L'}{L}(\lambda, \lambda,1/2+it)=-\frac12\log t+O(1)\qquad(T\to\infty).
\end{align}
\end{Lemma}
\proof
By the functional equation (\ref{Lerchfunc}), we have 
\begin{align*}
L(\lambda,\lambda,s)=&(2\pi)^{s-1}\Gamma(1-s)e^{\pi i\frac{1-s}{2}-2\pi i\lambda^2}\overline{L(\lambda,\lambda,1-\overline{s})}
\\&
\times\biggr(1
+\frac{e^{-\pi i(1-s)+2\pi i \lambda}L(\lambda,1-\{\lambda\},1-s)}
{\overline{L(\lambda,\lambda,1-\overline{s})}}\biggr).
\end{align*}
The logarithmic derivative gives
\begin{align}\label{logder}
\frac{L'}{L}(\lambda,\lambda,s)=\log 2\pi-\frac{\Gamma'}{\Gamma}(1-s)-\frac{\pi i}2-\overline{\frac{L'}{L}(\lambda,\lambda,1-\overline{s})}+E(\lambda,s),
\end{align}
where
\begin{align*}
E(\lambda,s)=\frac{\biggr(\frac{e^{-\pi i(1-s)+2\pi i \lambda}L(\lambda,1-\{\lambda\},1-s)}
{L(1-\lambda,\lambda,1-s)}\biggr)_s'}{1
+\frac{e^{-\pi i(1-s)+2\pi i \lambda}L(\lambda,1-\{\lambda\},1-s)}
{L(1-\lambda,\lambda,1-s)}}.
\end{align*}
For $0<\lambda,\alpha\le1$, we know that $L(\lambda, \alpha, 1-s)\ne0$, if $\sigma<-1$, moreover,  $L(\lambda, \alpha, 1-s)$ and its derivative have absolutely convergent Dirichlet series, if $\sigma<0$. Thus
$$E(\lambda,\sigma_2+it)\ll e^{-\pi t}\qquad(t\to\infty).$$
 By Stirling's formula, we get that
\begin{align}\label{gmm}
  \frac{\Gamma'}{\Gamma}(s) = \log s
  + O\left(|s|^{-1}\right)\quad( \Re(s) \geq 0, \ |s|\to\infty).
\end{align}
This proves the formula \eqref{left}.

We turn to the second part of Lemma \ref{lefthandside}. The expression (\ref{logder}) together with the formula \eqref{gmm} gives
\begin{align}\label{2re}
2\Re \frac{L'}{L}(\lambda,\lambda,1/2+it)=-\log t+\Re E(\lambda, 1/2+it)+O(1)\qquad(T\to\infty).
\end{align}
Next we consider the growth of $E(\lambda,1/2+it)$.  For $0<\lambda,\alpha\le1$, by Lemma~3 in \cite{gt} and by Cauchy's integral formula for the derivative, there is $B>0$ such that
 \begin{align}\label{eq:Fbound}
   L(\lambda, \alpha,1/2-iT)=O(T^B)\quad\text{and}\quad L'(\lambda, \alpha,1/2-iT)=O(T^B).
 \end{align}
 In view of the conditions of the lemma and the asymptotic formula \eqref{zeronumber} we have that  the distance $d$ from $1/2+it$ to the nearest zero of  $L(\lambda,\lambda, s)$ satisfies the inequalities 
 $$\exp\left(-A T/\log T\right)<d\ll1.$$
 Then Lemma \ref{1/L} yields
\begin{align*}
E(\lambda,1/2+it)\ll T^B\exp((-\pi+4CA)T+\log(3T)) .
\end{align*}
This finishes the proof of Lemma \ref{lefthandside}.

\endproof

 The following proposition will be important in the proof of Theorem \ref{dertest}. 

\begin{Proposition}\label{crho}
Let $0<\lambda\le1$, $T<t<T+U$, and $0<U\le T$. Let $A>0$ be such that $4AC<\pi$, where the constant $C$ is from Lemma \ref{1/L}. Let $\rho'$ be a zero of $L(\lambda,\lambda,s)$ such that $|\Re\rho'-1/2|<\exp(-AT/\log T)$ and  $T<\Im \rho'<T+U$. Assume that there are $0<\varepsilon<1$ and $\delta>0$ such that the function $L(\lambda,\lambda,s)$
has less than
\begin{align}\label{1log2}
\left [\frac{\varepsilon }{\log (2+\delta)} \log T\right]
\end{align}
zeros in the disc 
$|s-\rho'|\le\exp(-AT^{1-\varepsilon}/\log T)$.  Then, for sufficiently large $T$, there is a radius $r$,
\begin{align}\label{radiusr}
\exp\left(-\frac{A(1+\delta/3)T}{(2+\delta)\log T}\right)\le r\le\exp\left(-\frac{A(1+\delta/3)T^{1-\varepsilon}}{\log T}\right),
\end{align}
 such that $L(\lambda, \lambda, s)\ne0$ in the ring
\begin{align}\label{freering}
r^{(2+\delta)/(1+\delta/3)}\le |s-\rho'|\le r^{1/(1+\delta/3)}
\end{align}
and, for $|s-\rho'|=r$, $\sigma\le1/2$,
 \begin{align}\label{Nm}
\Re \frac{L'}{L}(\lambda,\lambda,s)\le -\frac12\log T+O(1)
\end{align}
\end{Proposition}
\proof
Let $r_k=\exp\left(-AT/((2+\delta)^k\log T)\right)$, $k=0,1,\dots,[\frac{\varepsilon }{\log (2+\delta)} \log T]$. 
By the condition \eqref{1log2} and Dirichlet's box principle, there is $j\in\{1,2,\dots,[\frac{\varepsilon }{\log (2+\delta)} \log T]\}$ such that the ring
\begin{align}\label{ring}
r_{j-1}<|s-\rho'|\le r_j
\end{align}
has no zeros of $L(\lambda,\lambda, s)$. Note that $r_j^{2+\delta}=r_{j-1}$. The function
$$
f(s)=\frac{L'}{L}(\lambda,\lambda,s)-\sum_{\rho\, :\, |\rho-\rho'|\le r_j^{2+\delta}}\frac{1}{s-\rho}
$$
is analytic in the disc $|s-\rho'|\le r_j$ and in this disc it has the Taylor expansion
\begin{align}\label{taylor}
f(s)=\sum_{n=0}^\infty a_n(s-\rho')^n.
\end{align}

We bound the coefficients $a_n$. Cauchy's integral formula for the derivative yields
\begin{align}\label{coeff}
a_n=\frac1{2\pi i}\int\limits_{|s-\rho'|= r_j}\frac{f(s)ds}{(s-\rho')^{n+1}}.
\end{align}
Lemma~3 in \cite{gt} gives that, for any $\sigma_0$, there is a positive constant $B$ such that 
 $L(\lambda, \lambda,s)=O(t^B)$ if $\sigma\ge\sigma_0$.  By the proof of Theorem 1 in \cite{gar} we see that, for any $t$ the modulus $|L(\lambda, \lambda,3+it)|$ is greater than some positive absolute constant. Therefore
by Lemma $\alpha$ from Titchmarsh \cite[\S 3.9]{Titchmarsh1986} and by the formula (\ref{zeronumber})  we obtain, for $|s-\rho'|\le r_j$,
\begin{align*}
\frac{L'}{L}(\lambda,\lambda,s)=\sum_{\rho\, :\, |\rho-(3+i\gamma')|\le 6}\frac{1}{s-\rho}+O(\log T).
\end{align*}
Then, in view of the definition of $f(s)$ using the zero free region (\ref{ring}), it follows that
$$
f(s)=\sum_{\rho\, :\, |\rho-(3+i\gamma')|\le 6\ \text{and}\atop |\rho-\rho'|> r_j}\frac{1}{s-\rho}+O(\log T).
$$
We apply the last expression to the formula \eqref{coeff}. For $n\ge1$ and $|\rho-\rho'|\ge r_j$, we have 
$$
\int_{|s-\rho'|=r_j} \frac{d s}{(s-\rho')^{n+1}(s-\rho)}=0
$$
and thus
\begin{align}\label{coeffexpr}
a_n
\ll
\int\limits_{|s-\rho'|= r_j}\frac{O(\log T)ds}{(s-\rho')^{n+1}}
\ll
r_j^{-n}\log T.
\end{align}

Now we choose $r=r_j^{1+\delta/3}$. Then expressions \eqref{taylor} and \eqref{coeffexpr} yield, for $|s-\rho'|=r$,
\begin{align*}
f(s)-a_0
=%&
\sum_{n=1}^\infty a_n(s-\rho')^n\ll \log T \sum_{n=1}^\infty r_j^{n\delta /3}
%\\ 
\ll%&
r_j^{\delta/3} \log T.
\end{align*}
By this we have, for $|s-\rho'|=r$,
\begin{align*}
\frac{L'}{L}(\lambda,\lambda,s)=%&
a_0+\sum_{\rho\, :\, |\rho-\rho'|\le r_j^{2+\delta}}\frac{1}{s-\rho}
%\\&
+O\left(r_j^{\delta/3} \log T\right).
\end{align*}
 Taking real parts we obtain, for $|s-\rho'|=r$, 
\begin{align}\label{takingrealparts}
\Re \frac{L'}{L}(\lambda,\lambda,s)
=%&
\Re a_0+\sum_{\rho\, :\, |\rho-\rho'|\le r_j^{2+\delta}}\frac{\sigma-\beta}{|s-\rho|^2}
%\\&
+O\left(r_j^{\delta/3}  \log T\right).%\nonumber
\end{align}

We will get an asymptotic formula for $\Re a_0$. We consider the sum over zeros in the formula \eqref{takingrealparts}. By inequalities $|\rho-\rho'|\le r_j^{2+\delta}$ and $|\Re\rho'-1/2|<r_0$ we see that, for $|s-\rho'|=r$, $1/2-(|\Re\rho'-1/2|+ r_j^{2+\delta})\le\sigma\le1/2$, and large $T$,
\begin{align}\label{numerator}
|\sigma-\beta|\le |\sigma-1/2|+|\Re \rho'-1/2|+|\Re \rho'-\beta|\le 4r_j^{2+\delta}
\end{align}
and
\begin{align}\label{denominator}
|s-\rho|^2\ge (|s-\rho'|-|\rho-\rho'|)^2=(r_j^{1+\delta/3}-r_j^{2+\delta})^2>r_j^{2+2\delta/3}/2.
\end{align}
The asymptotic formula \eqref{zeronumber} for the number of nontrivial zeros gives that there are $\ll\log T$ zeros $\rho$ such that $|\rho-\rho'|\le r_j^{2+\delta}$. Thus, for $|s-\rho'|=r$ and $1/2-(|\Re\rho'-1/2|+ r_j^{2+\delta})\le\sigma\le1/2$, we get
\begin{align*}
\sum_{\rho\, :\, |\rho-\rho'|\le r_j^{2+\delta}}\frac{\sigma-\beta}{|s-\rho|^2}\ll r_j^{\delta/3}\log T
\end{align*}
and
\begin{align}\label{isgasdino}
\Re \frac{L'}{L}(\lambda,\lambda,s)
=\Re a_0 +O\left(r_j^{\delta/3} \log T\right).
\end{align}
By \eqref{ring} we have that the ring $\{z : r_j^{2+\delta}<|z-\rho'|\le r_j\}$  has no zeros. Recall that $|s-\rho'|=r=r_j^{1+\delta/3}$. In view of this the distance  from $s=1/2+it$ to the nearest zero is 
$$\ge \min(r_j-r_j^{1+\delta/3}, r_j^{1+\delta/3}-r_j^{2+\delta})>r_0=\exp\left(-A T/\log T\right).
$$
Then the equality (\ref{logr11/2}) together with \eqref{isgasdino} gives
\begin{align}\label{rea}
\Re a_0=-\frac12 \log T+O(1).
\end{align}
By expressions \eqref{isgasdino} and \eqref{rea} we obtain that,  for $|s-\rho'|=r$ and $1/2-(|\Re\rho'-1/2|+ r_j^{2+\delta})\le\sigma\le1/2$, 
\begin{align}\label{final1}
\Re \frac{L'}{L}(\lambda,\lambda,s)
=-\frac12 \log T +O\left(1\right).
\end{align}

 If  $|s-\rho'|=r$ and $\sigma<1/2-(|\Re\rho'-1/2|+ r_j^{2+\delta})$, then we have 
 \begin{align*}
\sum_{\rho\, :\, |\rho-\rho'|\le r_j^{2+\delta}}\frac{\sigma-\beta}{|s-\rho|^2}\le 0
\end{align*}
and, in view of formulas \eqref{takingrealparts}, \eqref{rea},
\begin{align}\label{final2}
\Re \frac{L'}{L}(\lambda,\lambda,s)
\le-\frac12 \log T +O\left(1\right).
\end{align}
The expressions \eqref{final1} and \eqref{final2} together with the zero free region \eqref{ring} prove Proposition~\ref{crho}.

\endproof

 \proof[Proof of Theorem \ref{dertest}] Let
$$
R=\left\{s\in \mathbb C : T<t<T+U, -2<\sigma<\frac12
\right\}.
$$
We have that all the nontrivial zeros of $L(\lambda,\lambda, s)$ and $L'(\lambda,\lambda, s)$ lie to the right-hand side of the line $\sigma=-2$.
To start with, the idea is to consider the change of the argument of $L'/L(\lambda,\lambda, s)$ around the boundary of the region $R$. However, a problem occurs if $1/2+it$ is near to a zero of $L(\lambda,\lambda, s)$. Next, our goal is to exclude the zeros $\rho$ for which 
\begin{align}\label{almostoncritl}
|\beta-1/2|<\exp(-AT/\log T)\quad\text{and}\quad T<\gamma<T+U
\end{align}
 from the region $R$     using certain arcs which lie to the left-hand side of the line $\sigma=1/2$. We will use  Proposition \ref{crho}.

In this proof we always assume that the zero $\rho$ satisfies the inequalities \eqref{almostoncritl}. By the condition \eqref{Dcondition} of Theorem \ref{dertest}  there is $\delta>0$ such that the function $L(\lambda,\lambda,s)$
has less than
\begin{align*}
\left[\frac{\varepsilon }{\log (2+\delta)} \log T\right]
\end{align*}
zeros in the disc 
$|s-\rho|\le\exp(-AT^{1-\varepsilon}/\log T)$. 
Then, in view of Proposition \ref{crho}, for each such zero $\rho$, we define the disc
\begin{align*}
D(\rho,r)=\{s : |s-\rho|\le r\},
\end{align*}
where the radius $r$ is from Proposition \ref{crho}. Thus, for all $s$ such that $|s-\rho|=r$ and $\Re s\le1/2$, we have
\begin{align}\label{oncircle}
\Re \frac{L'}{L}(\lambda,\lambda,s)\le -\frac12\log T+O(1).
\end{align}

Let $S$ be the union of all the discs $D(\rho, r)$, where $\rho$ satisfies the inequalities \eqref{almostoncritl}. Note that each disc $D(\rho, r)$ from the set $S$ has a nonempty intersection with the critical line $\sigma=1/2$. For $T\le y\le T+U$ and $1/2+iy\not\in S$, formulas \eqref{logr11/2} and \eqref{radiusr} yield
\begin{align}\label{oncritline}
\Re \frac{L'}{L}(\lambda,\lambda,1/2+iy)\le -\frac12\log T+O(1),
\end{align}
provided that the area
\begin{align*}
V=&\left\{s : \left|s-(1/2+iT)\right|\le\exp(-AT/\log T), \ t\le T\right\}
\\
&\bigcup\left\{s : \left|s-(1/2+iT+iU)\right|\le\exp(-AT/\log T), \ t\ge T+U\right\}
\end{align*}
has no zeros of $L(\lambda, \lambda, s)$.

The vertical strips $|t-T|\le1$ and $|t-(T+U)|\le1$ contain $\ll \log T$ zeros of $L(\lambda, \lambda, s)$. Therefore, without loss of generality, we assume that $L(\lambda, \lambda, \sigma+iT)\ne0$, 
$L'(\lambda, \lambda,\sigma+iT)\ne0$ for $-2\le\sigma\le1/2$,  and $L(\lambda, \lambda, s)\ne0$ for $s\in V$. Further, we consider the
change of $\arg L'/ L(\lambda, \lambda,s)$ along the appropriately indented boundary
$R'$ of the region $R$.  More precisely, the upper, left, and lower sides
of $R'$ coincide with the upper, left, and lower boundaries of $R$. To
obtain the right-hand side of the contour $R'$, we take a curve $\psi$ defined as the boundary of the set $S\cup \{s : \sigma\ge1/2\}$, where this boundary is restricted to the strip $T\le t\le T+U$.

To prove the theorem, we will show that the change of $\arg  L'/ L(\lambda, \lambda,s)$ along the  contour $R'$ is $\ll\log T$. Let $R'_1$, $R'_2$, $R'_3$ and $R'_4$ denote the right, upper, left, and lower sides of the contour $R'$ accordingly.

We start from $\arg_{R'_1} L'/ L(\lambda, \lambda,s)$, where  $\arg_{R'_1} L'/ L(\lambda, \lambda,s)$ denotes the change of argument of $ L'/ L(\lambda, \lambda,s)$ along the right-hand side $R'_1$ of the contour $R'$. Formulas \eqref{oncircle} and \eqref{oncritline} give that
$$\left|\arg_{R'_1}\frac{L'}{L}(\lambda,\lambda,s)\right|<\pi.$$
Similarly,  the equality (\ref{left}) from Lemma \ref{lefthandside} gives
$$\left|\arg_{R'_3}\frac{L'}{L}(\lambda,\lambda,s)\right|<\pi.$$

Next we turn to horizontal sides $R'_2$ and $R'_4$.
 By standard arguments  using Jensen's
 theorem together with the bounds (\ref{eq:Fbound})  it is possible to show that (cf. \cite[inequality (7) and below]{gs2002} or Titchmarsh
 \cite[Section 9.4]{Titchmarsh1986})
$\arg_{R'_2}L(\lambda,\lambda,s)\ll\log T$, $\arg_{R'_2}L'(\lambda,\lambda,s)\ll\log T$, $\arg_{R'_4}L(\lambda,\lambda,s)\ll\log T$, and $\arg_{R'_4}L'(\lambda,\lambda,s)\ll\log T$.
This finishes the proof of Theorem \ref{dertest}.

\endproof

\section{Ending notes}\label{notes}

Here we discuss the curve $\psi$ from Theorem \ref{dertest}. Let $T<t\le T+U$.  
In the proof of  Theorem \ref{dertest},  we construct the curve $\psi$ which lies in the strip 
$$1/2-\exp(-AT^{1-\varepsilon}/\log T)\le \sigma \le 1/2.$$
Moreover, $\psi$ is constructed in a such way that the zero $\rho$ of  $L(\lambda,\lambda,s)$, lying in
\begin{align}\label{beta-}
1/2-\exp(-AT/\log T)\le \sigma\le 1/2,
\end{align} 
must also lie between the curve $\psi$ and the critical line $\sigma=1/2$. We expect that the location of the curve $\psi$ is not accidental and reflects interesting properties of the zeros of the Lerch zeta-function.
 In \cite{gt} we proved that if $\rho$ is a nontrivial zero of $L(\lambda,\lambda,s)$, then there is a radius $\exp(-A\gamma/\log \gamma)\le r\le \exp(-A\gamma/\log \gamma)\log^2\gamma$ such that the discs  
\begin{align}\label{discs}
|s-\rho|<r\quad\text{and}\quad |s-(1-\overline{\rho})|<r
\end{align}
 contain the same number of zeros. On the other hand, the calculations in \cite[Section 2]{gt} suggest that if $0<\lambda<1$, $\lambda\ne1/2$, and $\rho$ is a nontrivial zero of $L(\lambda,\lambda, s)$, then the symmetry described by the formula \eqref{discs} is not strict, namely, $1-\overline{\rho}$ is not a zero of $L(\lambda,\lambda, s)$. Moreover, if the discs in the expression \eqref{discs} intersect, then both discs possibly contain the same zero(s). From this we expect that the nontrivial zeros of $L(\lambda,\lambda, s)$, for $0<\lambda<1$, $\lambda\ne1/2$, can be classified into two classes, heuristically described as follows. One class contains zeros which are relatively far from the critical line. These zeros appear in almost symmetric pairs according to \eqref{discs}. Another class consists of zeros which are relatively near the critical line. They are almost symmetric to themselves (in view of \eqref{discs}). We expect that  the curve $\psi$ (or the  appropriate version of this curve lying ``nearest" to the critical line) from Theorem \ref{dertest} separates these two classes of zeros located in the left-hand side of the critical line.
 
Note that in Theorem \ref{dertest}, for $\lambda=1/2,1$, the curve $\psi$ can be constructed in at least two ways. One way is as described in the proof of Theorem \ref{dertest}. Another way is to choose $\psi(\tau)=1/2+iT+i\tau U$, where $\tau\in[0,1]$ (see Levinson and Montgomery \cite{Levinson1974} for the Riemann zeta function $L(1,1,s)$ and Y\i ld\i r\i m \cite{yildirim96} for the Dirichlet $L$-function $L(1/2,1/2,s)$). Clearly, if the generalized Riemann hypothesis is true and $\lambda=1/2,1,$ then in Theorem \ref{dertest} we can choose any curve $\psi$ which is located to the left from the critical line.

\end{document}